\documentclass[12pt]{amsart}
\usepackage{amscd,amssymb}
\usepackage[arrow,matrix]{xy}
\usepackage[plainpages,backref,urlcolor=blue]{hyperref}

\theoremstyle{plain}
\newtheorem{thm}[subsection]{Theorem}

\newtheorem{prop}[subsection]{Proposition}
\newtheorem{cor}[subsection]{Corollary}

\theoremstyle{definition}
\newtheorem{rk}[subsection]{Remark}

\newtheorem{ex}[subsection]{Example}

\numberwithin{equation}{section}
\newcommand{\OO}{{\mathcal O}}

\newcommand{\I}{{\mathcal I}}

\newcommand{\M}{{\mathcal M}}

\newcommand{\A}{{\mathcal A}}
\newcommand{\HH}{{\mathcal H}}

\newcommand{\B}{{\mathcal B}}
\newcommand{\CC}{{\mathcal C}}

\newcommand{\C}{\mathbb{C}}

\newcommand{\PP}{\mathbb{P}}

\DeclareMathOperator{\mult}{mult}

\begin{document}

\title [Unexpected curves and minimal degree of Jacobian relations]
{Unexpected curves in $\PP^2$, line arrangements, and minimal degree of Jacobian relations}

\author[A. Dimca]{Alexandru Dimca}
\address{Universit\'e C\^ ote d'Azur, CNRS, LJAD}
\email{dimca@unice.fr}



\begin{abstract} We reformulate a fundamental result due to Cook, Harbourne, Migliore and Nagel on the existence and irreduciblity of unexpected plane curves of a set of points $Z$ in $\mathbb{P}^2$, using the minimal degree of a Jacobian syzygy of the defining equation for the dual line arrangement $\mathcal A_Z$. Several applications of this new approach are given. In particular, we show that the irreducible unexpected quintics may occur only when the set $Z$ has the cardinality equal to  11 or 12, and describe five cases  where this happens.
\end{abstract}

\maketitle

\section{Introduction}
Let $Z=\{p_1, p_2, \ldots ,p_d\}$ be a finite set of $d$ points in the complex projective plane $\PP^2$.
One says that $Z$ admits unexpected curves of degree $j\geq 2$ if
$$h^0(\PP^2, \OO_{\PP^2}(j)\otimes \I(Z+(j-1)q))>\max\left(0, h^0(\PP^2, \OO_{\PP^2}(j)\otimes \I(Z))-{j \choose 2}\right),$$
where $q$ is a generic point in $\PP^2$, the fat point scheme $kq$ is defined by the $k$-th power of the corresponding maximal ideal sheaf $\I(q)$, and hence $\I(Z+(j-1)q)$ is the ideal sheaf of functions vanishing on $Z$ and vanishing of order $(j-1)$ at $q$, see \cite{CHMN,DMO,Just}. There is a more general definition, see \cite{HMNT,Just}, but in this note we consider only the special case described above.
Let $\A_Z:f_Z=0$ be the associated line arrangement in $\PP^2$ as in \cite{CHMN,DMO}. Let $(a_Z,b_Z)$ be the generic splitting type of the derivation bundle $E_Z$ associated to $\A_Z$, and let $m(\A_Z)$ be the maximal multiplicity of an intersection point in $\A_Z$. It is well known that $a_Z+b_Z=d-1$. For $i=1,2, \ldots,d$, let $Z_i=Z\setminus \{p_i\}$ be the set of $d-1$ points obtained from $Z$ be forgetting the point $p_i$, and let $\A_{Z_i}:f_{Z_i}=0$, $(a_{Z_i},b_{Z_i})$ and  $m(\A_{Z_i})$ be the corresponding objects associated with the set $Z_i$ as above.
With this notation,
the following fundamental result was established in \cite[Theorem 1.2]{CHMN},  \cite[Lemma 3.5 (a)]{CHMN}, \cite[Corollary 5.5]{CHMN} and \cite[Corollary 5.17]{CHMN},   see also \cite{DMO} for a discussion.
\begin{thm}
\label{thmA}
The set of points $Z$ admits an unexpected curve if and only if
$$m(\A_Z) \leq a_Z+1 <{d \over 2}.$$
If these conditions are fulfilled, then $Z$ admits an unexpected curve of degree $j$  if and only if
$$a_Z<j \leq d-a_Z-2.$$
The unexpected curve $C_q$ of minimal degree $j=a_Z+1$ and having a point of multiplicity $a_Z$ at a generic point $q$ is unique. Moreover $C_q$ is irreducible if and only if $a_Z=a_{Z_i}$ for all $i=1,2 \ldots,d$.
\end{thm}
 For larger values of $j$, the corresponding unexpected curves of degree $j$ are obtained from $\CC_q$ by adding
$j-a_Z-1$ lines passing through $q$, see \cite[Corollary 5.5]{CHMN}.
The curve $\CC_q$ itself, if not irreducible, is the union of some lines through $q$ and
an irreducible curve $\CC'_q$, having at $q$ a point of multiplicity $\deg(\CC_q')-1$.

Let $S=\C[x,y,z]$ be the polynomial ring in three variables $x,y,z$ with complex coefficients, and let $\A:f=0$ be an arrangement of  $d$ lines in the complex projective plane $\PP^2$. The minimal degree of a Jacobian syzygy for the polynomial $f$ is the integer $mdr(f)$
defined to be the smallest integer $r \geq 0$ such that there is a nontrivial relation
\begin{equation}
\label{rel_m}
 af_x+bf_y+cf_z=0
\end{equation}
among the partial derivatives $f_x, f_y$ and $f_z$ of $f$ with coefficients $a,b,c$ in $S_r$, the vector space of  homogeneous polynomials in $S$ of degree $r$.
The main result of this note is the following reformulation of Theorem \ref{thmA}.
\begin{thm}
\label{thmcorA}
The set of points $Z$ admits an unexpected curve if and only if
$$m(\A_Z) \leq mdr(f_Z)+1 <{d \over 2}.$$
If these conditions are fulfilled, then $Z$ admits an unexpected curve of degree $j$  if and only if
$$mdr(f_Z)<j \leq d-mdr(f_Z)-2.$$
Every unexpected curve $C_q$ of minimal degree $j=mdr(f_Z)+1$ is irreducible if and only if $mdr(f_Z)=mdr(f_{Z_i})$ for all $i=1,2 \ldots,d$.
\end{thm}
The advantage of having such a result comes from the wealth of information we have on the numerical invariant $mdr(f_Z)$, and on the relations between $mdr(f_Z)$ and $mdr(f_{Z_i})$ for various $i$, see
\cite{ADS,DMich,DSer}.
Using these results, we prove in this note some new results, and also give shorter proofs for some known results. In particular, the results about the irreducibility of the curves $C_q$ of minimal degree seem to be easier to prove using this new view-point.

In section 2 we recall some basic properties of the invariant $mdr(f)$,
and show in Proposition \ref{prop0} that a free line arrangement $\A:f=0$
is  supersolvable if and only if
 $mdr(f)=\min \{m-1,d-m\},$
where $m=m(\A)$ and $d=|\A|$.

In section 3, we show first that $a_Z=mdr(f_Z)$ when the set $Z$ admits unexpected curves, see Theorem \ref{aZ}, and use this equality to prove Theorem \ref{thmcorA} starting from Theorem \ref{thmA}. As an application, we give in Corollary \ref{corA} a short proof for the fact that the set of points $Z$ dual to the monomial arrangement $\A_m^0$, for $m\geq 5$, has irreducible unexpected curves of minimal degree $m+2$. This result was obtained first in
\cite[Proposition 6.12]{CHMN}. 

Then we prove in Proposition \ref{prop3C} a similar result for 
the set of points $Z$ dual to the full monomial arrangement $\M_m$, for $m \geq 4$. Note that the full monomial arrangement $\M_m$ is denoted by $\A_{3,m-2}^3$ in \cite{Just}, and the claim in Proposition \ref{prop3C} is
part of the claim in \cite[Theorem 6]{Just}. {\it However, the irreducible question does not seem to be  addressed in} \cite{Just}. By the results in \cite{HaHa}, the line arrangements $\M_m$ correspond exactly to the $m$-homogeneous supersolvable line arrangements having at least 3 modular points. We describe the unexpected curves for the 
$m$-homogeneous supersolvable line arrangements having 2 modular points in Proposition \ref{prop3D}.

We prove in Proposition \ref{prop2} that a set $Z$  with
$d=|Z|\leq 8$ never admits unexpected curves. In fact, if a set of $d=|Z| \leq 8$  admits unexpected curves, these curves must exist in degree $a_Z+1<{d \over 2} \leq 4$. However, as explained at the end of Example \ref{ex1}, such unexpected curves do not exist over a field of characteristic zero by the results in \cite{Ak,FGST}. Hence Proposition \ref{prop2} follows also from these results. We give below a short, independent proof of this claim. As a new result, we show that a set $Z$ in which at most 3 points are collinear does not admit  unexpected curves, see Proposition \ref{prop3}.
This fact fails in positive characteristics, see \cite[Example 2.4]{CHMN}.

In section 4, we reprove first the fact that a set $Z$  with
$d=|Z|=9$ admits  unexpected curves if and only if the associated line arrangement $\A_Z$ is projectively equivalent to the line arrangement $B_3=\M_4$. This result was proved first in \cite{FGST},  and we give a shorter proof in Proposition \ref{prop3B} using Theorem \ref{thmcorA}. 
Then we show in Theorem \ref{thm4A} that quintic irreducible unexpected curves may occur only when the set $Z$ consists of  11 or 12 points, and give 4 situations where such quintics occur.

In the final section we discuss several situations where we can add a new point $p'$ to $Z$  such that the new set $Z'=Z \cup \{p'\}$ also admits unexpected curves. First we discuss the arrangements $\A^1_m$ and $\A^2_m$, which interpolate between the arrangements $\A^0_m$ and
$\M_{m+2}=\A^3_m$ discussed above. We prove in both cases that the
corresponding dual set $Z$ has irreducible unexpected curves of minimal degree, claims that occur in \cite[Theorem 6]{Just} without a proof of the irreducibility.
In Proposition \ref{prop4} we discuss what happens when we add a generic point $p'$ to $Z$, and in Proposition \ref{prop5} we discuss what happens when we add a generic point $p'$ situated on a line in $\PP^2$
which contains a maximal number of points in $Z$, namely $m(\A_Z)$ points. This gives examples of sets $Z$ having unexpected curves, without being duals of free line arrangements, see Example \ref{ex2}.


The author would like to thank Takuro Abe, Lucja Farnik, Brian Harbourne, Lukas K\"uhne, Tomasz Szemberg and Justina Szpond for very useful discussions.

\section{Preliminaries}

Let $S=\C[x,y,z]$ be the polynomial ring in three variables $x,y,z$ with complex coefficients, and let $\A:f=0$ be an arrangement of  $d$ lines in the complex projective plane $\PP^2$. We denote by $n_j=n_j(\A)$ the number of intersection points in $\A$ of multiplicity $j$.  It is known that  $mdr(f)=0$ if and only if $n_d=1$, hence  $\A$ is a pencil of $d$ lines passing through one point. Moreover, 
$ mdr(f)= 1$ if and only if $n_d=0$ and $n_{d-1}=1$, hence  $\A$ is a 
near pencil, see for instance \cite{DIM}. For the definition and the basic properties of free and supersolvable line arrangements we refer to \cite{DHA}.

 Let $AR(f) \subset S^3$ be the graded $S$-module such, for any integer $j$, the corresponding homogeneous component $AR(f)_j$ consists of all the triples $\rho=(a,b,c)\in S^3_j$ satisfying \eqref{rel_m}. 
Let $\alpha$ be the minimum of the Arnold exponents $\alpha_p$ (alias singularity indices or log canonical thresholds, see Theorem 9.5 in \cite{Ko})  of the singular points $p$ of $\A$. The germ $(\A,p)$ is weighted homogeneous of type $(w_1,w_2;1)$ with $w_1=w_2={1 \over m_p}$, where $m_p$ is the multiplicity of $\A$ at $p$.
It is known that
\begin{equation}
\label{Ae}
\alpha_p={w_1}+{w_2}={2\over m_p},
\end{equation}
see for instance \cite[Formula (2.4.7)]{DS0}. 
With this notation, \cite[Theorem 9]{DS} can be restated in our setting as follows, see also \cite[Theorem 2.1]{DSer}.

\begin{thm}
\label{van}
Let $\A:f=0$ be an arrangement of $d$ lines in $\PP^2$ and $m=m(\A)$ be maximal multiplicity of an intersection point in $\A$.
Then $AR(f)_k=0$ for all 
$$k <{2 \over m}d-2.$$
Equivalently, one has
$$mdr(f) \geq {2 \over m}d-2.$$
\end{thm}

\begin{rk}
\label{rk3}
Let $\A:f=0$ be a line arrangement, and $p=(1:0:0)$  an intersection  point
on $\A$ of maximal multiplicity, say $m=\mult(\A,p)=m(\A)$. To this situation, one can associate a primitive Jacobian syzygy as explained in
\cite[Section 2.2]{DMich}. We recall this construction here.
Let $g=0$ be the equation of the subarrangement of $\A$ formed by the $m$ lines in $\A$ passing through $p$ and note that $g_x=0$. Then we can write $f=gh$ for some polynomial $h \in S$. The syzygy constructed as explained there is primitive
and has degree $r_p=d-m$, more precisely it is given by
\begin{equation}
\label{eqA}
\rho_p=(a,b,c)=(xh_x-d\cdot h,yh_x,zh_x),
\end{equation}
where $h_x$ denotes the partial derivative of $h$ with respect to $x$.
As shown in \cite[Theorem 1.2]{DMich}, the following cases are possible for $r=mdr(f)$.

\medskip

\noindent {\bf Case A:} $r=r_p=d-m$, in other words the constructed syzygy has minimal degree. In this case the arrangement $\A$, if free,  is in fact supersolvable, see Proposition \ref{prop0} below.
If $\A=\A_Z$, to have unexpected curves in this case we need
$$m \leq d-m+1 <{d \over 2}.$$
These two inequalities cannot both hold, so in this case there are no
unexpected curves.

\medskip 

\noindent {\bf Case B:} $r<r_p=d-m$, in other words the constructed syzygy has not minimal degree. Then the following two situations are possible.
\medskip 

\noindent {\bf Subcase B1:} $r=m-1$, and then $2m<d+1$ and $\A$ is free
with generic splitting type $a=m-1<b=d-m$. This case occurs exactly for the supersolvable line arrangements. Indeed, if $\A$ is supersolvable, with $m=m(\A)$ satisfying $2m \leq d+1$, then $m-1 \leq d-m$, and hence $r=m-1$, see \cite[Equation (2.2)]{AD2}. Conversely, a free line $\A:f=0$ arrangement such that $mdr(f)=m(\A)-1=m-1$ is supersolvable, see Proposition \ref{prop0} below.
Note that unexpected curves  occur in this case if and only if $d>2m$, see \cite[Theorem 3.7]{DMO}. As an example, the full monomial arrangement 
$$\M_m: f=xyz(x^{m-2}-y^{m-2})(y^{m-2}-z^{m-2})(z^{m-2}-x^{m-2})=0,$$
 is supersolvable, it has
$d=|\M_m|=3m-3$, $m(\M_m)=m$ and hence the condition $d>2m$ holds for any $m \geq 4$.

\medskip 

\noindent {\bf Subcase B2:} $m \leq r \leq d-m-1$, and then $2m <d$.
One example of this case is provided by the Fermat arrangements, a.k.a.
monomial arrangements 
$$\A^0_m: f_m=(x^m-y^m)(y^m-z^m)(z^m-x^m)=0,$$
see \cite{Just} for more information.
It is known that $m=m(\A^0_m)$, $d=3m$ and $r=mdr(f_m)=m+1$ for $m \geq 3$. The unexpected curves  occur in this case when $m \geq 5$ and are discussed in \cite[Proposition 6.12]{CHMN}. In particular, it is shown there that 
the unexpected curves of minimal degree $m+2$ are irreducible in this case. A new proof of this irreducibility is given below in Corollary \ref{corA}.
\end{rk}

\begin{prop}
\label{prop0}
A free line $\A:f=0$ arrangement  is supersolvable if and only if 
$$mdr(f)=\min \{m-1,d-m\},$$ 
where $m=m(\A)$ and $d=|\A|$. In particular, a line arrangement satisfying
$$m(\A) = mdr(f)+1 \leq {d \over 2}$$
is supersolvable.
\end{prop}
\proof
If $\A$ is supersolvable, then the claim follows from \cite[Proposition 3.2]{DStNSS}. Suppose now that $\A$ is free and
$$mdr(f)=\min \{m-1,d-m\}.$$
The next part of the proof was communicated to us by Takuro Abe, and uses
\cite[Proposition 4.2]{A}, where line arrangements in $\PP^2$ are regarded as central plane arrangements $\tilde A$ in $\C^3$.
Note that, since $\C^3 \setminus \tilde A=(\PP^2 \setminus \A) \times \C^*$,  one has
$$b_2(\C^3 \setminus \tilde A)= b_1(\PP^2 \setminus \A)+b_2(\PP^2 \setminus \A)=(d-1)+(m-1)(d-m).$$
The formula for $b_2(\PP^2 \setminus \A)$ used above is obtained using the formula for the Betti polynomial of a free arrangement $\A$ with exponents $(d_1,d_2)$, namely
$$B(\PP^2 \setminus \A;t)=(1+d_1t)(1+d_2t),$$
see for instance \cite[Theorem 8.3]{DHA}.
On the other hand, if we choose a flag 
$$X_3=\{0\} \subset X_2=L \subset X_1=P \subset X_0=\C^3,$$
where the line $L$ corresponds to a point $p$ in $\A$ of multiplicity $m$, and the plane $P$ corresponds to any line in $\A$ containing $p$, then 
$$\sum_{j=0}^2(|\tilde A_{X_{j+1}}|-|\tilde A_{X_{j}}|)|\tilde A_{X_{j}}|=0+(m-1)\cdot 1+(d-m)\cdot m=b_2(\C^3 \setminus \tilde A).$$
This equality implies, via \cite[Proposition 4.2]{A}, that the line arrangement $\A$ is supersolvable. The last claim follows, using the discussion in Remark \ref{rk3}, which shows that an arrangement satisfying $m-1=mdr(f)<d-m$ is free.
\endproof
{\it The interest of the supersolvable line arrangements in the study of unexpected curves} comes from the last claim in Proposition \ref{prop0}.
In addition, we have a very simple, purely combinatorial criterion for the existence of unexpected curves: a set $Z$, such that the dual line arrangement $\A_Z$ is supersolvable,
admits unexpected curves if and only 
\begin{equation}
\label{e1}
2m<d,
\end{equation}
where $d=|Z|=|\A_Z|$ and $m=m(\A_Z)$,  see \cite[Theorem 3.7]{DMO}.

We end this section with a side remark on irreducible curves in $\PP^2$, say of degree $d$ and having a point of multiplicity $d-1$.
\begin{prop}
\label{prop1}
Let $C:f=0$ be an irreducible curve of degree $d$ in $\PP^2$ having a singular point $p$ of multiplicity $d-1\geq 2$. Then the following hold.
\begin{enumerate}

\item $p$ is the only singular point of $C$;

\item $C$ is a rational curve;

\item the fundamental group $\pi_1(\PP^2 \setminus C)$ is abelian;

\item if the curve $C$ is cuspidal, i.e. if the singularity $(C,p)$ is irreducible, then $C$ is either free or nearly free.

\end{enumerate}
\end{prop}
\proof
The first two claims are well known. The third claim follows for instance from \cite[Corollary 4.3.8]{Dbook2}. The last claim follows from (3) using 
\cite[Corollary 3.2]{DStRIMS}.
\endproof

\begin{ex}
\label{ex1}
The $B_3$-arrangement  is a special case of the full monomial arrangement $\M_m$,
corresponding to $m=4$.
When $Z$ is the set of $9$ points dual to the $B_3$-arrangement, the curve $\CC_q$ is an irreducible quartic with an ordinary triple point at $q$. This was one of the motivating examples in developing this theory, and it has occured first in \cite{DIV}.
For more details, see Example 1.2 and Example 3.1  in \cite{DMO} as well as the detailed study in \cite{SzSz} where the explicit equations of the unexpected curves $C_q$ in this case are given.
{\it We do not know whether an unexpected curve can ever be cuspidal.}
Note that 9 is the minimal value for $|Z|$ such that $Z$ admits an unexpected curve, in view of Proposition \ref{prop2} below and 
for $|Z|=9$, the set $Z$ is unique up-to projective equivalence, see
\cite{FGST} and Proposition \ref{prop3B} below. However, if we work over algebraically closed fields of characteristic $p=2$, then there are unexpected cubic curves, which turn out to be cuspidal, see \cite[Example 2.4]{CHMN}, as well as \cite{Ak,FGST}. Such unexpected cubic curves do not exist over a field of characteristic zero by the results in \cite{Ak} and \cite[Theorem 1.4]{FGST}, see also Proposition \ref{prop3} below for an alternative proof.

\end{ex}

\section{The main results}
We have the following relation between the invariants $a_Z$ and $mdr(f_Z)$.
\begin{thm}
\label{aZ}
For any finite set $Z$, one has
$$a_Z \leq  \min \left (mdr(f_Z), \left \lfloor {d-1 \over 2} \right \rfloor \right).$$
Moreover, if $Z$ admits an unexpected curve, then
$$a_Z=mdr(f_Z).$$
\end{thm}
\proof
The first claim follows from  \cite[Proposition 3.2 (1)]{AD}. The second claim follows from Theorem \ref{thmA} which implies that
$a_Z <(d-2)/2$ in this case, and from \cite[Proposition 3.2 (2)]{AD}.
\endproof
\subsection{Proof of Theorem \ref{thmcorA}.}

If $Z$ admits an unexpected curve, then $a_Z=mdr(f_Z)$ and the claims, except the last one, are clear. On the other hand, if 
$$ mdr(f_Z)+1 <{d \over 2},$$
then it follows from \cite[Proposition 3.2 (2)]{AD} that $a_Z=mdr(f_Z)$,
and again the claims, except the last one, follow. For the last claim, note that $mdr(f_{Z_i}) \leq mdr(f_Z)$, see for instance \cite[Proposition 2.12]{ADS},
and hence 
$$mdr(f_{Z_i}) \leq mdr(f_Z)<{d \over 2} -1={(d-1)-1 \over 2}.$$
We claim that $a_{Z_i} = mdr(f_{Z_i})$ for any $i$, which would complete the proof.
Note that Theorem \ref{aZ} implies $a_{Z_i} \leq mdr(f_{Z_i})$, hence it is enough to show that the inequality $a_{Z_i} \leq  mdr(f_{Z_i})-1$ leads to a contradiction. Indeed, one has in this case
$$a_{Z_i} \leq mdr(f_{Z_i})-1 <{d\over 2} -2 ={(d-1)-3 \over 2}.$$
Using \cite[Proposition 3.2 (2)]{AD} for the line arrangement $\A_{Z_i}$,
we get $a_{Z_i} = mdr(f_{Z_i})$, hence a contradiction.
\bigskip

As a first application, we can give a shorter proof to the following known fact, see \cite[Proposition 6.12]{CHMN}.

\begin{cor}
\label{corA}
The set of points $Z_m$, dual to the monomial arrangement 
$$\A^0_m: f=(x^{m}-y^{m})(y^{m}-z^{m})(z^{m}-x^{m})=0,$$
admits unexpected curves of minimal degree $m+2$, for
$m\geq 5$, and all of them are irreducible.
\end{cor}
\proof
If $L$ is any line in $\A^0_m:f=0$, the number of intersection points on $L$
is exactly $m+1$. Recall that $d=|\A^0_m|=3m$ and $mdr(f)=m+1$.
We apply now \cite[Proposition 2.12]{ADS} to determine $mdr(f_L)$,
where $\A_L:f_L=0$ is the line arrangement obtained from $\A^0_m$ by deleting the line $L$. Since
$$|\A^0_m|-(m+1)=2m-1>m+1=mdr(f)$$
for $m\geq 3$, it follows that $mdr(f)=mdr(f_L)$. Theorem \ref{thmcorA} implies that the minimal degree unexpected curves are irreducible.
\endproof
One has also the following result, already stated in \cite[Theorem 6]{Just}.
\begin{prop}
\label{prop3C}
The  set of points $Z_m$, dual to the full monomial arrangement 
$$\M_m: f=xyz(x^{m-2}-y^{m-2})(y^{m-2}-z^{m-2})(z^{m-2}-x^{m-2})=0,$$
admits  unexpected curves of minimal degree $m$, for any 
$m\geq 4$, and all of them are irreducible.
\end{prop}
\proof
We know that an unexpected curve for $\M_m$ has degree $\geq mdr(f)+1=m$, and that $m \geq 4$ is a necessary and sufficient condition for the existence of such curves, see the discussion in Remark \ref{rk3}, Subcase B1. It remains to prove that such curves are irreducible, using
Theorem \ref{thmcorA}. Note that if we remove any line from $\M_m$,
the resulting arrangement $\A$ is still supersolvable, with $d=|\A|=3m-4$ and $m=m(\A)$. Since $m-1 \leq d-m=2m-4$ for $m \geq 3$, it follows that $mdr(f_{Z_i})=mdr(f)=m-1$. This completes the proof.
\endproof
By the results in \cite{HaHa}, the line arrangements $\M_m$ correspond exactly to the $m$-homogeneous supersolvable line arrangements having at least 3 modular points. The $m$-homogeneous supersolvable line arrangements having 2 modular points have been classified in \cite{AD2}, and are essentially subarrangements $\M_m^k(w)$ of $\M_m$, for $k=0,1, \ldots, m-3$, described as follows. In short, the line arrangement $\M_m^k(w)$ is obtained from $\M_m$ by deleting a number of lines passing through one fixed modular point, but not through the other two modular points.
Set $n=m-2$ and denote by $\mu_n$ the multiplicative group of the $n$-th roots of unity. For $1 \leq k < n$, let 
$$W(n,k)=\mu_n^k \setminus \Delta,$$
where 
$$\Delta=\{(w_1,\ldots,w_k)  \ : \ w_i \in \mu_n \text{ and } w_{j_1}=w_{j_2} \text{ for some } j_1 \ne j_2\}.$$
For $w=(w_1,\ldots,w_k) \in W(n,k)$ we define the line arrangement
$$\M_m^k(w): f(w)=xyz(x^n-y^n)(x^n-z^n)\prod_{j=1,k}(z-w_jy)=0$$
in $\PP^2$ for $1 \leq k <n$. We also set $\M_m^0(w):f(w)=xyz(x^n-y^n)(x^n-z^n)=0$, where $w$ is the empty sequence in this case.

\begin{prop}
\label{prop3D}
The set of points $Z_m^k(w)$, dual to the  line arrangement $\M_m^k(w)$,
admits unexpected curves of minimal degree $m$, for
any $m\geq 5$ and $2 \leq k \leq m-3$, and all of them are irreducible.
\end{prop}
\proof
First we note that $d=|\M_m^k(w)|=2m-1+k$ and $m=m(\M_m^k(w))$,
hence the inequality \eqref{e1} tells us that the  arrangement $\M_m^k(w)$
admits  unexpected curves if and only if $2 \leq k$. Consider from now the case $2 \leq k <m-2$, which implies $m \geq 5$.
The arrangement $\M_m^k(w):f(w)=0$ is supersolvable, and hence
$mdr(f(w))=m-1$, since $m-1 \leq d-m=m-1+k$ and we use Proposition \ref{prop0}. It follows that the unexpected curves have minimal degree $m$. To check that these curves are irreducible, we use Theorem \ref{thmcorA}. Let $L$ be a line in $\M_m^k(w)$ and denote by
$\M_L:f_L=0$ the line arrangement obtained from $\M_m^k(w)$ by deleting the line $L$.  If the line $L$ is not the line $L_x:x=0$, then the arrangement $\M_L$ is still supersolvable with at least a modular point of multiplicity $m_L=m$ and has $d_L=d-1=2m-2+k$ lines.
Hence 
$$mdr(f_L)=m_L-1=m-1 \leq m-2+k=d_L-m_L.$$
If $L=L_x$, the arrangement $\M_L$ is no longer supersolvable.
To determine $mdr(f_L)$, we use \cite[Proposition 2.12 (2)]{ADS}, since one has
$$|\M_m^k(w)|-|I_L|=(2m-1+k)-(k+2)=2m-3>m-1=mdr(f(w)),$$
where $I_L$ is the set of intersection points of $\M_m^k(w)$ on the line $L$. It follows that in both case one has $mdr(f_L)=mdr(f(w))$, and by 
Theorem \ref{thmcorA} the unexpected curves of minimal degree are irreducible.
\endproof

The following two results say that, if a set $Z$ admits unexpected curves, then the associated line arrangement $\A_Z$ has to be rather complicated. Both results use Theorem \ref{van}, and hence apply only in characteristic zero.

\begin{prop}
\label{prop2}
A set of points $Z$ with $d=|Z| \leq 8$ does not admit unexpected curves.
\end{prop}
\proof
 We prove only the case $d=8$, since the other cases are easier and can be treated in a completely similar way. Assume that $Z$ has unexpected curves.
Using Theorem \ref{thmA} we get 
$$m(\A_Z) \leq a_Z+1 <4$$
and hence $a_Z \leq 2$ and $m(\A_Z) \leq 3$. Using Theorem \ref{van}
we get that 
$$mdr(f_Z) \geq {2\over 3} 8-2={10 \over 3}.$$ Hence $mdr(f_Z) >3$.
On the other hand we know that $a_Z=mdr(f_Z)$  by Theorem \ref{aZ}.
This contradiction proves our claim.
\endproof

\begin{prop}
\label{prop3}
A set of points $Z$ such that at most 3 points in $Z$ are collinear does not admit unexpected curves. In other words, a set of points $Z$ such $m(\A_Z) \leq 3$, does not admit unexpected curves. In particular, there are no unexpected cubics.
\end{prop}
Note that the case $m(\A_Z) = 2$ was treated in \cite[Corollary 6.8]{CHMN}, and a new, quick proof for this result can also be obtained using exactly the same argument as below.
\proof
It is enough, by the above remark, to treat the case $m(\A_Z) = 3$.
Then Theorem \ref{van} implies
$$mdr(f_Z) \geq {2 \over 3}d-2.$$
If $Z$ admits unexpected curves, we have in addition by Theorem \ref{thmcorA}.
$$ mdr(f_Z)+1 <{d \over 2}.$$
But
$${2 \over 3}d-2 <{d \over 2}-1$$
holds only for $d \leq 5$, and in this range $Z$ does not admit unexpected curves by Proposition \ref{prop2}. The existence of unexpective cubics would imply $m(\A_Z) \leq 3$ by Theorem \ref{thmA}, which is not possible.
\endproof

\section{On quartic and quintic unexpected curves in $\PP^2$}

We have seen in Proposition \ref{prop3} that there are no unexpected cubics.
The following result says that the unexpected quartics may occur only in a unique situation. It was first proved in \cite{FGST}, but our proof seems shorter.
\begin{prop}
\label{prop3B}
A set of points $Z$ admits unexpected quartics  only if $d=|Z|=9$.
Moreover, a set of points $Z$ with $d=|Z|=9$  admits unexpected curves if and only if $Z$ is projectively equivalent to the set of 9 points dual to the $B_3$-arrangement described in Example \ref{ex1} above, and then the unexpected curves are irreducible quartics.
\end{prop}
\proof
If $Z$ admits unexpected quartics, we get, using Theorem \ref{thmA} and Proposition \ref{prop3}, the following
$$4 \leq m(\A_Z) \leq mrd(f_Z)+1 =4<{d \over 2}.$$
Hence $mdr(f_Z) =3$, $m(\A_Z)=4$ and $d \geq 9$. Proposition \ref{prop0} implies that the arrangement $\A_Z$ is supersolvable, and then \cite[Theorem 1.1]{AD2}  implies
$$d=|\A_Z| \leq 3m(\A_Z)-3=9.$$
Hence the only possibility is $d=9$.
The numbers $n_k$, of the intersection points of multiplicity $k$ in a line arrangement $\A$ with $|\A|=d$, satisfy a number of relations. The easiest of them is  the following.
	\begin{equation}
\label{eqSum}
\sum_{k\geq 2} n_k{k \choose 2}={d \choose 2}, 
\end{equation}
where $d=|\A|$.
For a line arrangement $\A:f=0$, one has
\begin{equation}
\label{eqSS1}
\tau(\A)\leq (d-1)^2-r(d-r-1),
\end{equation} 
where
$$\tau(\A)=\sum_{k\geq 2}n_k(k-1)^2$$
and 
$r=mdr(f)$, and equality holds in \eqref{eqSS1} if and only if $\A$ is free, see \cite{Dmax,dPW}. In our case, $d=9$ and $r=3$, so we get two equations
$$ n_2+3n_3+6n_4=36 \text{ and } n_2+4n_3+9n_4=49.$$
The only solutions of these two equations, consisting of non-negative integers $n_j$, are the following four  vectors
$$(n_2,n_3,n_4) \in \{(9,1,4),(6,4,3),(3,7,2), (0,10,1)\}.$$

A highly non-trivial restriction on these numbers is given by the Hirzebruch inequality, valid for non trivial line arrangements (i.e. for line arrangements not a pencil or a near pencil), see \cite{H}:
\begin{equation}
\label{eqHir}
n_2 +\frac{3}{4}n_3 -d\geq  \sum_{k >4}(k-4)n_k.
\end{equation}
Using \eqref{eqHir}, it follows that the vector $(n_2,n_3,n_4)$ of our
arrangement of 9 lines with unexpected curves can be only
$(9,1,4)$ and $(6,4,3)$. Using the classification of supersolvable arrangements with at least 3 modular points given in \cite{HaHa}, our claim is proved.
\endproof

\begin{rk}
\label{rk3B}
The last two triples $(n_2,n_3,n_4)$ in the proof above can also be discarded using the fact that, for a supersolvable line arrangement $\A$, one has
$$n_2(\A) \geq {|\A| \over 2},$$
see \cite{A3,T1}.
\end{rk}

\begin{cor}
\label{cor3E}
A set of points $Z$ with $d=|Z|=10$ does not admit unexpected curves.
\end{cor}
\proof
Assume that $Z $ admits unexpected curves.
Using Proposition \ref{prop3} we see that $m(\A_Z) \geq 4$. On the other hand, Theorem \ref{thmcorA} implies
$$m(\A_Z) \leq mdr(f_Z)+1 <{d \over 2}=5.$$
It follows that $m(\A_Z)=4$ and $mdr(f_Z)=3$, and hence the unexpected curves a quadrics. This contradicts Proposition \ref{prop3B}, and proves our claim.

\endproof

Now we show that the unexpected quintics can occur only in a very limited number of situations. The arrangement $\A_m^2$, for $m\geq 3$,  is obtained from the full monomial arrangement $\M_{m+2}$ by deleting a line joining two modular points, see Proposition \ref{prop4} below for details.

\begin{thm}
\label{thm4A}
A set of points $Z$ in the plane $\PP^2$ admits unexpected irreducible quintics only if $d=|Z|\in \{11,12\}$.
Moreover, for each value of $d=|Z|\in \{11,12\}$, there are at least the following two possibilities.
\begin{enumerate}
\item For $d=12$, the set $Z$ can be the dual of the line arrangement
$\M_5$, or of the Hessian line arrangement given by
$$\HH: f_{\HH}=xyz\left( (x^3+y^3+z^3)^3-27x^3y^3z^3\right)=0.$$
\item For $d=11$, the set $Z$ can be the dual of the line arrangements
$\M_5^2$ and $\A_3^2$, obtained from $\M_5$ by deleting a line, or of the arrangement $\HH'$ obtained from of the Hessian  arrangement $\HH$ by deleting a line.
\end{enumerate}
\end{thm}
The set $Z$, dual to the Hessian line arrangement $\HH$, can be described explicitly as follows: $Z$ consists of the three points
$ (1:0:0)$, $(0:1:0)$ and $(0:0:1)$ and of the nine points
$$(w_1:w_2:1),$$
where $(w_1,w_2) \in \mu_3^2$, see \cite[Exercise 5.5]{DHA}.
\proof
If $Z$ admits unexpected quintics, we get as above, using Theorem \ref{thmA} and Proposition \ref{prop3},
$$4 \leq m(\A_Z) \leq mrd(f_Z)+1= 5<{d \over 2}.$$
This implies $mdr(f)=4$ and $d\geq 11$. There are two cases to consider.

Assume first that $m(\A_Z)=4$. Then Theorem \ref{van} implies that
$$4 \geq  {d \over 2} -2,$$
and hence $d \leq 12$. Assume next that $m(\A_Z)=5$, and apply Proposition \ref{prop0} to get that $\A_Z$ is supersolvable. Then apply
\cite[Theorem 1.1]{AD2}  which yields
$$d=|\A_Z| \leq 3m(\A_Z)-3=12.$$

The fact that the supersolvable line arrangements $\M_5$ and $\M_5^2$ admit irreducible unexpected quintic curves follows from Proposition
\ref{prop3C} and Proposition
\ref{prop3D}. The claim for the supersolvable line arrangement $\A_3^2$
follows from  Proposition \ref{prop4} below.

The Hessian line arrangement $\HH$ and the line arrangement $\HH'$ are both free with $mdr(f)=4$. The claim for $\HH$ is well known, see for instance \cite[Example 8.6 (iii)]{DHA}. The claim for $\HH'$ follows for instance using Terao's addition-deletion theorem, in the form stated in \cite[Theorem 3.2 (2)]{ADS}. Indeed, any line $L$ in $\HH$ has exactly 2 double points and 3 intersection points of multiplicity 4 on it, see for instance \cite[Exercise 5.5]{DHA}. 

This discussion implies the irreducibility claim for the unexpected curves in the case of the
Hessian line arrangement $\HH$. For $\HH'$,
we apply now \cite[Proposition 2.12 (2)]{ADS} to determine $mdr(f_L)$,
where $\HH'_L:f_L=0$ is the line arrangement obtained from $\HH'$ by deleting a line $L \in \HH'$. As said above, any line
$L \in \HH'$ has at most 5 intersection points on it, hence
\cite[Proposition 2.12 (2)]{ADS} can be applied to this situation to get
$$mdr(f_L)=mdr(f_{\HH'})=4.$$
\endproof

\begin{rk}
\label{rk4A}
We do not know whether the above 5 cases listed in Theorem \ref{thm4A} are the only ones where unexpected irreducible quintics may occur. Using the database of line arrangements with characteristic polynomials splitting over the integers, given in \cite{Lukas}, we know that we have listed all the cases for $m(\A_Z)=5$, as well as all the cases for $m(\A_Z)=4$ and $\A_Z$ a free arrangement.
When $m(\A_Z)=4$, $|\A_Z|=11$ and $\A_Z$ is not a free arrangement, we can show that $\A_Z:f_Z=0$ is a plus-one generated line arrangement using \cite[Theorem 2.3]{DSt3syz}. Moreover, since $\A_Z$ is not free, it follows from \eqref{eqSS1} that
$$\tau(\A_Z)< 10^2-4 \times 6=76.$$
On the other hand, \cite[Theorem 1.1 (ii)]{Dmin} implies that
$$\tau(\A_Z)\geq 10\times 6 +6+6+1=73.$$
Hence, when $d=11$, any non free line arrangement $\A_Z$ with $m(\A_Z)=mdr(f_Z)=4$ has to satisfy 
$$\tau(\A_Z) \in \{73,74,75\}.$$
Conversely, let $\A:f=0$ be a line arrangement such that $d=|\A|=11$, $m(\A)=4$, $\A$ not free, and $\tau(\A) \in \{73,74,75\}.$ 
{\it Lukas K\"uhne has informed me that there are 333 combinatorially distinct types of line arrangements satisfying these conditions: 296 types with $\tau(\A)=73$, 30 types with $\tau(\A)=74$ and 7 types with $\tau(\A)=75$.
This result was obtain using the database of all matroids of size 11 given in \cite{MaMo}. The 7 types with $\tau(\A)=75$ correspond to free arrangements $\A$, hence give nothing new for our problem.
 It remains in principle to study the possible values of the analytic invariant $mdr(f)$ and to list all the cases with $mdr(f)=4$ when
$\tau(\A) \in \{73,74\}.$ In all the examples tested so far, Lukas K\"uhne has found $mdr(f) >4$.}

Using a stronger form of the inequality \eqref{eqSS1}, namely
\begin{equation}
\label{eqSS1bis}
\tau(\A) \leq(d-1)(d-r-1)+r^2-{ 2r-d+2 \choose 2},
\end{equation}
valid for $2r\geq d$, see \cite{dPW},
we see that $r=mdr(f)$ has to satisfy $4 \leq r \leq 6$.
When $r=6$, we get $\tau(\A)\leq 73$, and hence by our assumptions
$\tau(\A)=73$. 
\end{rk}

\section{Adding a new point to $Z$}
First we revisit some results stated in \cite[Theorem 6]{Just}. Starting with the monomial arrangement $\A^0_m:f^0=0$, denoted $\A^0_3(m)$ in \cite{Just},
one can add the line $L_x:x=0$ and get the new line arrangement
$$\A^1_m: f^1=xf^0=x(x^m-y^m)(y^m-z^m)(z^m-x^m)=0.$$
\begin{prop}
\label{propF1}
The line arrangement $\A^1_m$ is free, but not supersolvable, with exponents $(m+1,2m-1)$
and the corresponding dual set of points $Z$ admits unexpected curves
of degree $j$ for $m\geq 4$ any integer $j$ satisfying
$$m+2\leq j\leq 2m-2.$$
The unexpected curves of minimal degree $j=m+2$ are all irreducible.
\end{prop}
\proof
First we apply \cite[Proposition 2.12]{ADS} for the line arrangement $\A'=\A^0_m$ and $H=L_x$. Since the set $I_H$ of intersection points of 
$\A^1_m$ on $H$ has cardinal $m+2$, and since
$$|\A'|-|I_H|=3m-(m+2)=2m-2>m+1=mdr(f^0)$$
for $m \geq 4$, it follows that $mdr(f^1)=m+1$.
Using this equality, it is easy to check that $\A^1_m$ is free using the equation \eqref{eqSS1}. The other claims, except the irreducibility claim,
follow from Theorem \ref{thmcorA}. Finally we address the irreducibility question.
A line $L$ in $\A^1_m$ has either $m+2$ intersection points
if $L=L_x$, or just $m+1$ intersection points when $L \ne L_x$.
We apply now \cite[Proposition 2.12]{ADS} to determine $mdr(f_L)$,
where $\A_L:f_L=0$ is the line arrangement obtained from $\A^1_m$ by deleting the line $L$. Since
$$|\A^1_m|-|I_L|\geq (3m+1)-(m+2)=2m-1>m+1=mdr(f^1)$$
for $m\geq 3$, where $I_L$ denotes the set of intersection points of $\A^1_m$ situated on the line $L$.
 It follows that $mdr(f^1)=mdr(f_L)$. Theorem \ref{thmcorA} implies that the minimal degree unexpected curves are irreducible.
\endproof
Starting now with the monomial arrangement $\A^1_m:f^1=0$, 
one can add the line $L_y:y=0$ and get the new line arrangement
$$\A^2_m: f^2=yf^1=xy(x^m-y^m)(y^m-z^m)(z^m-x^m)=0.$$
\begin{prop}
\label{propF2}
The line arrangement $\A^2_m$ in $\PP^2$ is supersolvable and  has a unique modular point, namely $p=(0:0:1)$. In particular, $\A^2_m$ is free with exponents $(m+1,2m)$
and the corresponding dual set of points $Z$ admits unexpected curves
of degree $j$ for $m \geq 3$ and any integer $j$ satisfying
$$m+2\leq j\leq 2m-1.$$
The unexpected curves of minimal degree $j=m+2$ are all irreducible.
\end{prop}
\proof Since $\A^2_m$ is clearly supersolvable, with $(0:0:1)$ as modular point, all the claims except the claim about irreducibility are proved using Theorem \ref{thmcorA}.
As above, denote by $\A_L:f_L=0$ is the line arrangement obtained from $\A^2_m$ by deleting the line $L$. It is enough to show that $mdf(f^2)=mdr(f_L)$. If $I_L$ denotes the set of intersection points of $\A^2_m$ situated on the line $L$, it is clear that $|I_L|\leq m+2$.
Since
$$|\A^2_m|-|I_L|\geq (3m+2)-(m+2)=2m>m+1=mdr(f^2)$$
for $m\geq 2$, the result follows by \cite[Proposition 2.12]{ADS}.
\endproof

In the following two results, we add a point $p'$ to $Z$, and hence the corresponding dual line $L'$ to the arrangement $\A_Z$. In both cases, the unexpected curves of minimal degree are not irreducible, as follows using Theorem \ref{thmcorA}.
\begin{prop}
\label{prop4}
Assume that the set of points $Z$ satisfies the stronger condition
$$m(\A_Z) \leq mdr(f_Z)+{3 \over 2} <{d \over 2}.$$
Let $p'$ be a generic point in $\PP^2$ and consider the new set $Z'=Z \cup \{p'\}$. Then $Z'$ admits an unexpected curve of degree $j$,  for any integer $j$ such that 
$$mdr(f_Z)+1<j \leq d-mdr(f_Z)-2.$$
\end{prop}
\proof
The point $p'$ gives by duality a generic line $L'$. Hence the arrangement $\A_{Z'}$ is given by adding a generic $L'$ to $\A_Z$.
Using \cite[Proposition 4.11]{ADS}, it follows that
$$m(\A_{Z'})=m(\A_Z), \ mdr(f_{Z'})=mdr(f_Z)+1 \text{ and } |\A_{Z'}|=|\A_Z|+1.$$
The claim follows using Corollary \ref{corA}.
\endproof

The point $p'$ in Proposition \ref{prop4} is generic if and only if $p'$ is not situated on any line $\overline {p_ip_j}$ determined by two distinct points $p_i,p_j \in Z$. Note also that the multiplicity $m=m(\A_Z)$ is exactly the maximal number of points in $Z$ which are collinear.
Let $p_{i_1},p_{i_2}, \ldots, p_{i_m}$ be a maximal set of collinear points in $Z$ and let $L$ be the line determined by these points. With this notation, we have the following result.

\begin{prop}
\label{prop5}
Assume that the set of points $Z$ satisfies the stronger condition
$$m(\A_Z) \leq mdr(f_Z)+{3 \over 2} <{d \over 2}.$$
Let $p'$ be a generic point on the line $L$ defined above and consider the new set $Z'=Z \cup \{p'\}$. Then $Z'$ admits an unexpected curve of degree $j$,  for any integer $j$ such that 
$$mdr(f_Z)+1<j \leq d-mdr(f_Z)-2.$$
\end{prop}
\proof
The point $p'$ gives by duality a line $L'$, which is generic in the pencil of lines passing through the common intersection point $p_L$ of the lines
$L_j$, dual to the points $p_{i_j}$, for $j=1,\ldots,m$. In fact, $p_L$ is the point dual to the line $L$.
The arrangement $\A_{Z'}$ is given by adding the line $L'$ to $\A_Z$.
Using \cite[Proposition 4.10]{ADS}, it follows that
$$m(\A_{Z'})=m(\A_Z)+1, \ mdr(f_{Z'})=mdr(f_Z)+1 \text{ and } |\A_{Z'}|=|\A_Z|+1.$$
Indeed, the case (3) in \cite[Proposition 4.10]{ADS} cannot occur, as explained in Remark \ref{rk3}, Case A.
The claim follows using Corollary \ref{corA}.
\endproof

\begin{ex}
\label{ex2}
When $Z$ is the set of $3m$ points dual to the Fermat $\A^0_m$-arrangement considered in Remark \ref{rk3}, Subcase B2, the conditions in Proposition \ref{prop4} are fulfilled for any
$m \geq 6$. Note that the arrangement $\A_{Z'}=\A_m \cup L'$ from
Proposition \ref{prop4}
 is far from being a free arrangement. Indeed, the global Tjurina number of
the arrangement $\A_{Z'}=\A_m \cup L'$ is given 
$$\tau(\A_{Z'})=\tau(\A_m)+3m=7m^2-3m+3.$$
On the other hand, the global Tjurina number of
a free arrangement $\B_m$ of $3m+1$ lines with $mdr(\B_m)=m+2$ is given by the formula \eqref{eqSS1} and hence
$$\tau(\B_m)=9m^2-(m+2)(2m-2)=7m^2-2m+4 >\tau(\A_{Z'}).$$
Hence $\A_{Z'}=\A_m \cup L'$ gives rise to countable many examples of sets $Z'$ admitting unexpected curves, and such that 
the corresponding arrangements  $\A_{Z'}$ are not free, and in particular not supersolvable.
\end{ex}


\begin{thebibliography}{00}

\bibitem{A}  T. Abe, Restrictions of free arrangements and the division theorem, in: "Perspectives in Lie Theory", Springer INdAM Series 19 (2017), 389--401.

\bibitem{A3}
T. Abe, Double points of supersolvable and divisionally free line arrangements in the projective plane, arXiv:1911.10754.


\bibitem{AD}  T. Abe, A. Dimca, On the splitting types of bundles of logarithmic vector fields along plane curves, Internat. J. Math. 29 (2018), no. 8, 1850055, 20 pp.

\bibitem{AD2} T. Abe, A. Dimca, On complex supersolvable line arrangements, arXiv: 1907.12497.



\bibitem{ADS}  T. Abe, A. Dimca, G. Sticlaru, Addition-deletion results for the minimal degree of logarithmic derivations of arrangements, arXiv:1908.06885.

\bibitem{Ak} S. Akesseh, Ideal containments under flat extensions and interpolation on linear systems in $\PP^2$, PhD thesis, University of Nebraska-Lincoln (2017).

\bibitem{T1}  B. Anzis, S.O. Toh\u aneanu, On the geometry of real and complex supersolvable line arrangements, J. Combin. Theory, Ser. A 140(2016), 76--96.


\bibitem{Lukas} M. Barakat, R. Behrends, C. Jefferson, L. K\" uhne, M. Leuner, On the generation of rank 3 simple matroids with an application to Terao's freeness conjecture, 	arXiv:1907.01073.


\bibitem{SzSz}T. Bauer, G. Malara, T. Szemberg, J. Szpond, Quartic unexpected curves and surfaces,  Manuscripta Math. (2018). https://doi.org/10.1007/s00229-018-1091-3.


\bibitem{CHMN} D. Cook, B. Harbourne, J. Migliore, U. Nagel,  Line arrangements and configurations of points with an unexpected geometric property. {Compositio Math.} {154}(2018), 2150--2194.

\bibitem{DIV} R. Di Gennaro, G. Ilardi, J. Vallès,  Singular hypersurfaces characterizing the Lefschetz
properties. J. Lond. Math. Soc. (2) 89(2) (2014), 194--212. 

\bibitem{DMO} M. Di Marca, G. Malara, A. Oneto, Unexpected curves arising from special line arrangements, Journal of Algebraic Combinatorics, https://doi.org/10.1007/s10801-019-00871-0


\bibitem{Dbook2}  A. Dimca, Singularities and Topology of Hypersurfaces, Universitext, Springer Verlag, New York,
1992.



\bibitem{DHA}  A. Dimca, Hyperplane Arrangements: An Introduction, Universitext, Springer-Verlag, 2017.

\bibitem{DMich}  A. Dimca, Curve arrangements, pencils, and Jacobian syzygies,  Michigan Math. J. 66 (2017), 347--365.

\bibitem{Dmax}  A. Dimca, Freeness versus maximal global Tjurina number for plane curves, Math. Proc. Cambridge Phil. Soc.  163 (2017), 161--172.

\bibitem{Dmin}  A. Dimca, On the minimal value of global Tjurina numbers for line arrangements, European J. Math. 
https://doi.org/10.1007/s40879-019-00373-0 

\bibitem{DIM} A. Dimca, D. Ibadula, A. M\u acinic, Numerical invariants and moduli spaces for line arrangements, arXiv:1609.06551, Osaka J. Math. (to appear)

\bibitem{DS0}
A. Dimca, M. Saito:
Some remarks on limit mixed Hodge structure and spectrum, {An. \c St. Univ. Ovidius Constan\c ta} 22(2) (2014), 69-78.

\bibitem{DS} {A. Dimca,  M. Saito}: Generalization of theorems of Griffiths
and Steenbrink to hypersurfaces
with ordinary double points, Bull. Math. Soc. Sci. Math. Roumanie, 60(108) (2017), 351--371.


\bibitem{DSer} A. Dimca, E. Sernesi,  Syzygies and logarithmic vector fields along plane curves, Journal de l'\'Ecole polytechnique-Math\'ematiques 1(2014), 247-267.


\bibitem{DStRIMS}
A. Dimca, G. Sticlaru, 
Free and nearly free curves vs. rational cuspidal plane curves, {Publ. RIMS Kyoto Univ. } {54} (2018), 163--179.

\bibitem{DStNSS} A. Dimca, G. Sticlaru, On supersolvable and nearly supersolvable line arrangements,
Journal of Algebraic Combinatorics 50 (2019), 363-378.

\bibitem{DSt3syz} A. Dimca, G. Sticlaru, Plane curves with three syzygies, minimal Tjurina curves, and nearly cuspidal curves, 
Geometriae Dedicata 10.1007/s10711-019-00485-7




\bibitem{dPW} A.A. du Plessis,  C.T.C. Wall, Application of the theory of the discriminant to highly singular plane curves, {Math. Proc. Camb. Phil. Soc.},  {126} (1999), 259-266. 


\bibitem{FGST} \L. Farnik, F. Galuppi, L. Sodomaco, W. Trok, On the unique unexpected quartic in $\PP^2$,  Journal of Algebraic Combinatorics, DOI 10.1007/s10801-019-00922-6.

\bibitem{HaHa} K. Hanumanthu, B. Harbourne, Real and complex supersolvable line arrangements in the projective plane, arXiv: 1907.07712.

\bibitem{HMNT} B. Harbourne, J. Migliore, U. Nagel, and Z. Teitler. Unexpected hypersurfaces and where to
find them, arXiv:1805.10626, accepted for publication in Michigan Math. J.

\bibitem{H} F.~Hirzebruch,
{\em Arrangements of lines and algebraic surfaces}, Arithmetic and geometry, Vol. II, Progr.
Math. 36, Birkhauser, Boston, Mass., 1983, 113--140.

\bibitem{Ko} J. Koll\' ar: Singularities of pairs, Algebraic Geometry, Santa Cruz, 1995; {Proceedings
of Symposia in Pure Math.} vol. 62, AMS, 1997, pages 221--287.

\bibitem{MaMo} Matsumoto, Y., Moriyama, S., Imai, H. et al. Matroid Enumeration for Incidence Geometry. Discrete Comput Geom 47 (2012), 17-43.



\bibitem{Just}  J. Szpond, Fermat-type arrangements, arXiv:1909.04089.


\end{thebibliography}
\end{document}